\documentclass[a4paper,12pt]{article}

\pdfoutput=1

\usepackage[paper=letterpaper,margin=1in]{geometry}
\usepackage{amsmath,amssymb,amsfonts,epsfig,cite,setspace,bigstrut,longtable,array,breqn,color}
\usepackage{hyperref,caption,setspace}
\usepackage{amsthm}
\usepackage{pdfpages,lipsum}
\usepackage{todonotes}
\usepackage{cancel}

\usepackage[OT2,T1]{fontenc}
\DeclareSymbolFont{cyrletters}{OT2}{wncyr}{m}{n}
\DeclareMathSymbol{\Sha}{\mathalpha}{cyrletters}{"58}

\usepackage{mathtools}
\newcommand\myeq{\stackrel{\mathclap{\normalfont\mbox{?}}}{=}}

\newenvironment{red}{\relax\color{red}}{\relax}
\newenvironment{blue}{\relax\color{blue}}{\hspace*{.5ex}\relax}
\newcommand{\ber}{\begin{red}}
\newcommand{\er}{\end{red}}
\newcommand{\beb}{\begin{blue}}
\newcommand{\eb}{\end{blue}}

\newcommand{\sage}{{\sc SageMath }} 
\newcommand{\mth}{{\sc Mathematica }}
\begin{document}


\vskip 0.25in

\newcommand{\nn}{\nonumber}
\newcommand{\tr}{\mathop{\rm Tr}}
\newcommand{\comment}[1]{}
\newcommand{\cM}{{\cal M}}
\newcommand{\cW}{{\cal W}}
\newcommand{\cN}{{\cal N}}
\newcommand{\cH}{{\cal H}}
\newcommand{\cK}{{\cal K}}
\newcommand{\cZ}{{\cal Z}}
\newcommand{\cO}{{\cal O}}
\newcommand{\cA}{{\cal A}}
\newcommand{\cB}{{\cal B}}
\newcommand{\cC}{{\cal C}}
\newcommand{\cD}{{\cal D}}
\newcommand{\cT}{{\cal T}}
\newcommand{\cV}{{\cal V}}
\newcommand{\cE}{{\cal E}}
\newcommand{\cF}{{\cal F}}
\newcommand{\cX}{{\cal X}}
\newcommand{\IA}{\mathbb{A}}
\newcommand{\IP}{\mathbb{P}}
\newcommand{\IQ}{\mathbb{Q}}
\newcommand{\IH}{\mathbb{H}}
\newcommand{\IR}{\mathbb{R}}
\newcommand{\IC}{\mathbb{C}}
\newcommand{\IF}{\mathbb{F}}
\newcommand{\IV}{\mathbb{V}}
\newcommand{\II}{\mathbb{I}}
\newcommand{\IZ}{\mathbb{Z}}
\newcommand{\re}{{\rm~Re}}
\newcommand{\im}{{\rm~Im}}

\let\oldthebibliography=\thebibliography
\let\endoldthebibliography=\endthebibliography
\renewenvironment{thebibliography}[1]{%
\begin{oldthebibliography}{#1}%
\setlength{\parskip}{0ex}%
\setlength{\itemsep}{0ex}%
}%
{%
\end{oldthebibliography}%
}

\newtheorem{theorem}{\bf THEOREM}
\def\thetheorem{\thesection.\arabic{theorem}}
\newtheorem{proposition}{\bf PROPOSITION}
\def\thetheorem{\thesection.\arabic{proposition}}
\newtheorem{observation}{\bf OBSERVATION}
\def\thetheorem{\thesection.\arabic{observation}}
\newtheorem{conjecture}{\bf CONJECTURE} 
\def\thetheorem{\thesection.\arabic{CONJECTURE}}

\theoremstyle{definition}
\newtheorem{definition}{\bf DEFINITION} 
\def\thetheorem{\thesection.\arabic{DEFINITION}}
\newtheorem{example}{\bf EXAMPLE} 
\def\thetheorem{\thesection.\arabic{EXAMPLE}}
\newtheorem{remark}{\bf REMARK} 
\def\thetheorem{\thesection.\arabic{REMARK}}

\def\theequation{\thesection.\arabic{equation}}
\newcommand{\setall}{\setcounter{equation}{0}
        \setcounter{theorem}{0}}
\newcommand{\setequation}{\setcounter{equation}{0}}
\renewcommand{\thefootnote}{\fnsymbol{footnote}}

\newcommand{\seteq}{\mathbin{:=}}
\newcommand{\GL}{\operatorname{GL}}
\newcommand{\Sp}{\operatorname{Sp}}
\newcommand{\USp}{\operatorname{USp}}
\newcommand{\GSp}{\operatorname{GSp}}
\newcommand{\SU}{\operatorname{SU}}
\newcommand{\SO}{\operatorname{SO}}
\newcommand{\End}{\operatorname{End}}

\begin{titlepage}

~\\
\vskip 1cm

\begin{center}
{\Large \bf Machine-Learning Arithmetic Curves}
\end{center}
\medskip

\renewcommand{\arraystretch}{0.5} 

\vspace{.4cm}
\centerline{
{\large Yang-Hui He, Kyu-Hwan Lee, Thomas Oliver}
}
\vspace*{3.0ex}

\vspace{10mm}

\begin{abstract}
We show that standard machine-learning algorithms may be trained to predict certain invariants of low genus arithmetic curves. Using datasets of size around $10^5$, we demonstrate the utility of machine-learning in classification problems pertaining to the BSD invariants of an elliptic curve (including its rank and torsion subgroup), and the analogous invariants of a genus 2 curve. Our results show that a trained machine can efficiently classify curves according to these invariants with high accuracies ($>0.97$). For problems such as distinguishing between torsion orders, and the recognition of integral points, the accuracies can reach $0.998$.
\end{abstract}

\end{titlepage}

\begin{spacing}{1}
\tableofcontents
\end{spacing}

\section{Introduction}

In this article we build on recent work by the present authors, namely: \cite{HLOa}, \cite{HLOb}. In the latter, we presented experiments demonstrating the capacity of machine-learning to predict basic invariants of algebraic number fields. Many of the invariants studied there appear together in the analytic class number formula:
\begin{equation}\label{eq.ACNF}
\lim\limits_{s \to 1} (s-1) \zeta_F(s) = \frac{2^{r_1} (2\pi)^{r_2} \, \mathrm{Reg}_F \, h_F}{w_F\sqrt{|\Delta_F|}} \ ,
\end{equation}
in which $\zeta_F(s)$ is the Dedekind zeta function of a number field $F$, and, $(r_1,r_2)$ is the signature, $\mathrm{Reg}_F$ is the regulator, $h_F$ is the class number, $\Delta_F$ is the discriminant, and $w_F$ is the number of roots of unity in $\mathcal{O}_F$. 

Recall that the set $E(\mathbb{Q})$ of rational points on an elliptic curve $E$ defined over $\mathbb{Q}$ defines a finitely generated abelian group. We will denote the rank of this group by $r$, and the torsion subgroup by $E(\mathbb{Q})_{\mathrm{tors}}$. Associated to an elliptic curve $E$ over $\mathbb{Q}$, one has the $L$-function $L(E,s)$ which conjecturally vanishes to order $r$ at $s=1$ with the leading Taylor coefficient given by the following equation analogous to~\eqref{eq.ACNF}:
\begin{equation}\label{eq.BSDF}
\frac{L^{(r)}(E,1)}{r!}~\myeq~\frac{\#\Sha(E/\mathbb{Q})\Omega(E/\mathbb{Q})\mathrm{Reg}(E/\mathbb{Q})\prod_{p}c_p}{(\#E(\mathbb{Q})_{\mathrm{tors}})^2},
\end{equation}
in which the symbol $~\myeq~$ indicates that the equality is conjectural, $\Sha(E/\mathbb{Q})$ is the Tate--Shafarevich group, $\Omega(E/\mathbb{Q})$ is the global period, $\mathrm{Reg}(E/\mathbb{Q})$ is the regulator, $c_p$ is the Tamagawa number at a prime $p$ (cf. \cite[Appendix C, Conjecture~16.5]{silverman} for details). 
More generally, with appropriate modifications, it is possible to replace $E/\mathbb{Q}$ in equation~\eqref{eq.BSDF} by the Jacobian of a smooth, projective, geometrically integral curve defined over a global field. 

In this paper, we will apply machine-learning techniques to predict - with varying levels of success - the invariants appearing in equation~\eqref{eq.BSDF}.

Equation~\eqref{eq.BSDF} is of course the famous conjecture of Birch and Swinnerton-Dyer (BSD), which is a Millennium prize problem. Aside from the BSD conjecture, there are many interesting open questions regarding the numbers appearing in equation~\eqref{eq.BSDF}. For example, it is unknown whether or not the Tate--Shafarevich group appearing on the right-hand side is finite, cf. \cite[Chapter~X,~Conjecture~4.13]{silverman}, and there is no rigorous method known for its computation. Nevertheless, when the Tate--Shafarevich group is finite, its order is known to be a square \cite[Chapter~X,~Theorem~4.14]{silverman}. As for the left-hand side, it is not yet known whether or not, as $E$ varies, the set of ranks $r$ is bounded (a heuristic model suggesting that this set might be bounded was presented in \cite{PPVW}). Furthermore, in a suitable asymptotic sense, it is conjectured that 50\% of elliptic curves have rank 0, 50\% have rank 1, and $0\%$ have rank $\geq2$. \cite[Conjecture~B]{Goldfeld}, \cite{KSa,KSb}.

For a broad introduction to machine-learning, see \cite{ML,Hastie}. In this paper we will utilise logistic regression, naive Bayes, and random forest classifiers, which are reviewed in \cite[Sections~4.4,~6.6.3,~15]{Hastie}. Machine-learning algorithms require the training of classifiers using large sets of data, and we obtain our data sets from \cite{lmfdb}. Previous elliptic curve machine-learning experiments were documented in \cite{Alessandretti:2019jbs}, as part of a recent programme of machine-learning various structures in mathematics \cite{ml1,ml2,ml3}.
The key difference in the present article is that, whilst \cite{Alessandretti:2019jbs} utilised Weierstrass coefficients as training data (which had enormous variation in magnitude), we will here lean more heavily on the Euler factors of $L$-functions. We observe this to be much more successful, allowing even for extrapolation to elliptic curves with conductors in ranges beyond those in the training dataset. 

This paper also studies genus 2 curves over $\mathbb{Q}$. There are some important differences between elliptic curves and those of genus 2. For example, there are known to be only finitely many rational points on a genus 2 curve. In contrast, an elliptic curve with positive rank has infinitely many rational points. On the other hand, the rational points on the Jacobian of a genus 2 curve form a finitely generated abelian group which could be infinite. By the rank of a genus 2 curve, we mean the rank of its Jacobian. For elliptic curves, classification by rank is essentially binary because curves of higher rank do not provide enough training data (as per the conjecture mentioned above). On the other hand, there is a significant proportion of genus 2 curves with rank 2 - thus allowing for a ternary classification. 

We remark that there is some comparison to be made with \cite{HLOa}, in which we studied machine-learning of a particular classification problem arising from the Sato--Tate conjecture for hyperelliptic curves.
There, we found that naive Bayes classifiers could distinguish, with accuracies into $0.99 \sim 1.00$ range, the Sato--Tate groups of hyperelliptic curves with genus 1 or 2, using a small number of Euler-factors. The experimental results of this paper show that the same method is just as powerful for other invariants. 
It is interesting that a method as simple as naive Bayes could do so well for various invariants in number theory. It might be suggesting that mathematics is more workable with machine-learning than the real world where data sets could be dimmed or distorted by various noises. 

We conclude this introduction with an overview of what is to come. In Section~\ref{s:notation} notations are fixed with brief explanations of some concepts. In Section~\ref{s:strategy}, the generation of training data and the experimental set-up are explained. In Section~\ref{s:elliptic} we document the experimental outcomes for elliptic curves. In Section~\ref{s:genus2}, we do the same for genus 2 curves. Finally, in Section~\ref{s:outlook}, we offer some concluding remarks and tentative directions for further research.

\subsection*{Acknowledgements}

We thank \'Alvaro Lozano-Robledo, Andrew Sutherland and Chris Wuthrich for helpful discussions. YHH is indebted to STFC UK, for grant ST/J00037X/1, KHL is partially supported by a grant from the Simons Foundation (\#712100), and TO acknowledges support from the EPSRC through research grant EP/S032460/1.

\section{Notation}\label{s:notation}

We use the following notation throughout:

\begin{description}

\item[Elliptic curve] defined over $\mathbb{Q}$ is denoted by $E$. The set $E(\mathbb{Q})$ of rational points defines a finitely generated abelian group;

\item[Genus 2 curve] defined over $\mathbb{Q}$ is denoted by $C$. The curve $C$ is assumed to be smooth, projective, and geometrically integral. The set $C(\mathbb{Q})$ of rational points is finite;

\item[Jacobian] of $C$ is denoted by $J$. The Jacobian is a two-dimensional Abelian variety defined over $\mathbb{Q}$. The set $J(\mathbb{Q})$ of rational points defines a finitely generated Abelian group;

\item[Rank] of $E$ (resp. $C$) denoted by $r_E$ (resp. $r_C$) is the rank of the finitely generated abelian group $E(\mathbb{Q})$ (resp. $J(\mathbb{Q})$). If the rank is 0 (resp. positive), then there are finitely many (resp. infinitely many) rational points;

\item[Torsion subgroup] of $E(\mathbb{Q})$ (resp. $J(\mathbb{Q})$) is denoted by $E(\mathbb{Q})_{\mathrm{tors}}$ (resp. $J(\mathbb{Q})_{\mathrm{tors}}$);

\item[Cyclic Group] of order $n$ is denoted $\mathtt C_n$. The torsion subgroup of $E(\mathbb{Q})$ is a product of cyclic groups;

\item[Good primes] of a variety $X$ defined over $\mathbb{Q}$ are those primes $p\in\mathbb{Z}$ such that $X$ has an integral model whose reduction modulo $p$ defines a smooth variety of the same dimension. A good prime for $C$ is a good prime for $J$, but the converse is not necessarily true;

\item[Bad primes] of a variety $X$ defined over $\mathbb{Q}$ are those primes $p\in\mathbb{Z}$ which are not good. The bad reduction types of an elliptic curve are reviewed in \cite[Section~VII.5]{silverman};

\item[Conductor] of $E$ (resp. $C$) denoted by $Q_E$ (resp. $Q_C$) is a positive integer of the form $\prod p^{e_p}$ in which $p$ varies over the bad primes for $E$ (resp. $J$). The power $e_p$ to which a bad prime $p$ appears depends on the reduction type, cf. \cite[Section~IV.10]{silverman2};

\item[Tate--Shafarevich group] of $E$ (resp. $J$) denoted by $\Sha(E/\mathbb{Q})$ (resp. $\Sha(J/\mathbb{Q})$) is a torsion Abelian group and measures the extent to which the Hasse principle fails to hold, cf. \cite[Section~X.4]{silverman}.

\end{description}

\section{Methodology}\label{s:strategy}

In this section we explain our experimental set-up. In particular, we construct the training and validation sets from appropriate data and overview the machine-learning strategies used.

\subsection{Euler factors}\label{s:Ltraining}

Let $X$ be a smooth, projective, geometrically connected curve of genus $g\in\{1,2\}$. For each good prime $p$ of $X$, we define the local zeta function to be:
\begin{equation}
Z(X/\mathbb F_p; T) = \exp \left ( \sum_{k=1}^\infty \frac{\#X\left(\mathbb{F}_{p^k}\right)T^k}{k} \right ) . 
\end{equation}
It is well-known that the local zeta function can be written in the form
\begin{equation}\label{eq.genus2ZL}
Z(X/\mathbb F_p; T) = \frac{L_p(X,T)}{(1-T)(1-pT)} , 
\end{equation}
where $L_p(T) \in \mathbb Z[T]$ is a  polynomial of degree $2g$ with constant term $1$.

\begin{example}
If $X=E$ is an elliptic curve defined over $\mathbb{Q}$ and $p$ is a good prime for $E$, then:
\begin{equation}\label{eq.EllipticEuler}
L_p(E,T)=1 - a_p T + p T^2 , \ \ (p\text{ good}),
\end{equation}
where
\begin{equation}\label{eq.ap}
a_p=p+1-\#E\left(\mathbb{F}_{p}\right).
\end{equation}
For a bad prime $p$, we also define $a_p$ as in equation~\eqref{eq.ap}. Using \sage \cite{sage}, we may compute a large amount of $a_p$ quickly. For $i\in\mathbb{Z}_{>0}$, let $p_i$ denote the $i$th prime. For a positive integer $N$, we introduce the vector:
\begin{equation}\label{eq.Lvector}
v_L(E)=(a_{p_1},\dots,a_{p_N})\in\mathbb{Z}^{N}, \ \ N\in\mathbb{Z}_{\geq1}.
\end{equation}
In practice, we will take $N$ to be $100$, $200$, $300$ or $500$. We note that the $100$th prime is $541$, $200$th is $1223$, the $300$th is $1987$, and the $500$th is $3571$.
\end{example}
\begin{example}
If $X=C$ is a smooth projective geometrically connected genus $2$ curve defined over $\mathbb{Q}$ and $p$ is a good prime for $C$, then:
\begin{equation}\label{eq.Genus2Euler}
L_p(C,T)=1+a_{1,p}T+a_{2,p}T^2+a_{1,p}pT^3+p^2T^4, \ \ a_{1,p},a_{2,p}\in\mathbb{Z}.
\end{equation}
For a bad prime $p$, we will simply use the convention
\begin{equation}
(a_{1,p},a_{2,p})=(0,p).
\end{equation}
Using \sage \cite{sage}, we may compute  $(a_{1,p},a_{2,p})$. For a positive integer $N$, we introduce the vector:
\begin{equation}\label{eq.Lvector-1}
v_L(C)=\left((a_{1,p_2},a_{2,p_2}),\dots,(a_{1,p_{N+1}},a_{2,p_{N+1}})\right)\in\left(\mathbb{Z}^2\right)^{N}, \ \ N\in\mathbb{Z}_{\geq1},
\end{equation}
where we do not include $p_1=2$ as it is always bad.
In practice, we will take $N=200$.
\end{example}
Given a finite set $\mathcal{F}$ of $X$  and an invariant $I(X)$ for each $X$, we associate the following labeled dataset:
\begin{equation}\label{eq.Ldata}
\mathcal{D}_{L,X}=\{v_L(X)\rightarrow I(X):X\in\mathcal{F}\}.
\end{equation}
We will refer to the entries in $v_L(X)$ as the Euler coefficients of $X$. When $X$ is an elliptic curve (resp. a genus $2$ curve), the Euler coefficients are integers (resp. pairs of integers).

\subsection{Experimental strategy}\label{sec:strategy}

\begin{enumerate}

\item Let $\mathcal{F}$ be a finite set of elliptic curves (resp. smooth projective geometrically connected genus 2 curves). The choice of $\mathcal{F}$ depends on the experiment. For example, $\mathcal{F}$ could be the set of elliptic curves (resp. genus 2 curves) over $\mathbb{Q}$ with conductor less than some bound and rank in the set $\{0,1\}$.

\item For an elliptic curve $E$ (resp. genus 2 curve $C$) in $\mathcal{F}$, let $I(E)$ (resp. $I(C)$) denote an invariant of interest. For example, $I(E)$ (resp. $I(C)$) could be the rank of $E(\mathbb{Q})$ (resp. $J(\mathbb{Q})$).

\item Generate datasets of the form $\mathcal{D}=\{v(X)\rightarrow I(X):X\in\mathcal{F}\}$, where $\mathcal{D}$ is as in~\eqref{eq.Ldata}\footnote{In exceptional circumstances, we will in fact construct different datasets in place of $\mathcal{D}$. We will do this, for example, in the investigation of particularly accurate classifiers as in Section~\ref{sec:conjectures}, and in an attempt to improve on a poorly performing classifier as in Section~\ref{sec:ECTS}. Such a digression from convention will always be clearly indicated.}. We will take $N$ to be one of: $100, 200,300,500$. We stress at this point that $N$ is an absolute constant, and does not vary with the curves in $\mathcal{D}$.

\item Choose a subset $\mathcal{T}\subset\mathcal{D}$ and denote its complement by $\mathcal{V}=\mathcal{D}-\mathcal{T}$. We will refer to $\mathcal T$ as the training dataset, and $\mathcal V$ as the validation dataset. It is important that the training set and validation set have no intersection so as not to over-fit the machine-learning.
We will not typically specify $\mathcal T$, or its size relative to $\mathcal{D}$, as the choice will not impact significantly on the results. See also step 7.
 
\item Train a classifier on the set $\mathcal{T}$ with a standard supervised-learning algorithm. In this paper we will use naive Bayes, random forests, and logistic regression - see \cite[Sections~4.4,~6.6.3,~15]{Hastie}. We implement the algorithms using \mth \cite{wolfram}.

\item For all curves $X$ in $\mathcal{V}$, ask the classifier to determine $I(X)$. We record the precision and confidence, which together constitute a good measure of accuracy and performance of the machine. The precision and confidence are real numbers in the interval $[0,1]$, and the aspiration is that both are close to $1$. By precision, we mean the proportion of predictions in agreement with \cite{lmfdb}, the validity of which is discussed in \cite[Reliability~of~elliptic~curve~data~over~$\mathbb{Q}$,~Reliability~of~genus~2~curve~data~over~$\mathbb{Q}$]{lmfdb}. By confidence, we mean the Matthew's correlation coefficient \cite{matthews}. The confidence value is an extra check intended to minimize false positives and false negatives.

\item Repeat steps 4 to 6 for different choices of $\mathcal{T}$. The precision and confidence values recorded below are representative of several repetitions. 

\end{enumerate}

\section{Elliptic curves}\label{s:elliptic}

In this section we describe our experimental results for elliptic curves defined over $\mathbb{Q}$. For standard algorithms used in the computation of the invariants discussed below, the reader is referred to \cite{cremona}.  To perform the experiments in this section, we downloaded data from \cite[Elliptic~curves~over~$\mathbb{Q}$]{lmfdb}, the completeness of which is discussed in \cite[Completeness~of~elliptic~curve~data~over~$\mathbb{Q}$]{lmfdb}. We note that the Hasse--Weil $L$-function of an elliptic curve $E$ is an invariant of its isogeny class, and so we in fact downloaded a representative curve for each isogeny class. On the LMFDB, an isogeny class is represented by an optimal curve, and hence our data sets are generated from optimal curves only. In general, the torsion order, torsion structure and the number of integral points, considered in Sections \ref{sec:torsorder} - \ref{sec:integral},  are not uniquely determined by an isogeny class.   

\subsection{Rank}\label{sec:ECrank}

\begin{table}[h!!!]
\begin{center}
{\footnotesize \begin{tabular}{|c|c|c|c|c|c|}
\hline
$Q_E$ training range &$N$ & $Q_E$ validation range & |Data| = $\#\{E\}$&Precision&Confidence\\
\hline
$[1,1\times10^4]$  & 100 &$[1,1\times10^4]$  & 
                     $1.6\times10^4$ ($\times2$) & 0.977&0.955\\
\hline
" & 300 & " & " & 0.991&0.982\\
\hline
$[2\times10^4+1,3\times10^4]$&300&$[2\times10^4+1,3\times10^4]$& $1.7\times10^4$ ($\times2$) &0.964&0.922\\
\hline
"&500&"&"&0.971&0.941\\
\hline
$[1,1\times10^4]$ &300& $[2\times10^4+1,3\times10^4]$ &"&0.924&0.848\\
\hline
\end{tabular}
}
\end{center}
\caption{{\sf
The above table shows the precision and confidence of a logistic regression classifier when asked to distinguish elliptic curves over $\mathbb{Q}$ with rank $0$ from those with rank $1$. The classifier is trained on $E$ with conductor $Q_E$ in the ranges specified by the first column, using the number of Euler factors given in the second column. The classifier is verified on $E$ with conductor in the ranges specified by the third column.
}}
\label{t:rank}
\end{table}

Recall that we denote by $r_E$ the rank of an elliptic curve $E$.
 
It is conjectured that, in a rigorous sense, $50\%$ of elliptic curves over $\mathbb{Q}$ have rank $0$, $50\%$ have rank 1, and $0\%$ have higher rank, cf. \cite[Conjecture~B]{Goldfeld}, \cite{KSa,KSb}. 
Furthermore, it is known that if $r_E\leq1$ then $r_E$ is equal to the order of vanishing of $L(E,s)$ at $s=1$.
It is therefore expedient to consider this as a binary classification problem using the vectors $v_L$ defined by Euler factors as in \eqref{eq.Lvector}. For different ranges of conductor $Q_E$, we established a balanced dataset of size $\sim 2 \times10^4~(\times 2)$ for rank 0 and rank 1.

Trying several standard classifiers, we find that logistic regression worked best and the results are summarized in Table~\ref{t:rank}.
We see that the accuracies are in the high 0.90s, which is reassuring that a machine learns ranks of elliptic curves.
What is of particular interest is the last line in the table, where we trained on 300 Euler factors for conductors in the range from 1 to $10^4$ but validated on those in the range from $2\times10^4+1$ to $3\times10^4$, and still achieved a 0.92 precision.

The results show that the number of Euler factors needed for high precision is about $3 \sqrt{\max\{Q_E\}}$ in the range of $Q_E$ we considered. 
We also note that a logistic regression classifier also performed best in distinguishing the ranks of algebraic number fields \cite{HLOb} when number fields were presented through defining polynomials. On the other hand, when trained on Weierstrass coefficients as in \cite{Alessandretti:2019jbs}, no classifier was able to accurately predict the rank of an elliptic curve.

\subsection{Torsion order}\label{sec:torsorder}

\begin{table}[h!!!]
\begin{center}
{\small \begin{tabular}{|c|c|c|c|c|}
\hline
$Q_E$ range& |Data| $= \#\{E\}$ & $N$ &Precision&Confidence\\
\hline
$[1,3\times10^4]$ & $3.75\times10^4$ ($\times2$)& 500 &0.9997&0.9995\\
\hline
\end{tabular}
}
\end{center}
\caption{{\sf
The above table shows the precision and confidence of a naive Bayes classifier when asked to distinguish elliptic curves over $\mathbb{Q}$ with torsion order $1$ from those with torsion order $2$.  The classifier is trained on a random sample of curves with conductor in the range specified by the first column, and verified on those which remain.
}}
\label{t:torsorder}
\end{table}

The torsion group of an elliptic curve over $\mathbb{Q}$ has order at most $16$ \cite[Chapter~VII,~Theorem~7.5]{silverman}.  Currently, there are too few data-points on LMFDB to experiment with torsion order $>2$. 
Thus, we perform supervised machine-learning of the form $\{ v_L \} \to   ( |E(\mathbb{Q})_{\mathrm{tors}}| = 1 \mbox{ or } 2)$, whereby predicting the torsion group being trivial or not, using the Euler coefficients alone. To be clear, we do not restrict the rank of a curve in this experiment.
We established a balanced dataset of size $\sim 4 \times10^4~(\times 2)$ for torsion order 1 and 2 together.
A naive Bayes classifier was used and the results are summarized in Table~\ref{t:torsorder}.
We see that the accuracies are extremely good, using 500 Euler coefficients. We note that the naive Bayes classifier appeared also in \cite{HLOa}.  We will revisit this experiment in Section~\ref{sec:conjectures}.

\subsection{Torsion structure} \label{ss:torsion-structure}

\begin{table}[h!!!]
\begin{center}
{\small \begin{tabular}{|c|c|c|c|c|}
\hline
$Q_E$ range & |Data| =  $\#\{E\}$ & $N$ &Precision&Confidence\\
\hline
$[1,1\times10^6]$ &  $5.4\times10^3$ ($\times2$) & 500 &0.885&0.789\\
\hline
\end{tabular}
}
\end{center}
\caption{{\sf
The above table shows the precision and confidence of a random forest classifier when asked to distinguish elliptic curves over $\mathbb{Q}$ such that $E(\mathbb{Q})_{\mathrm{tors}}\cong \mathtt C_4$ from those such that $E(\mathbb{Q})_{\mathrm{tors}}\cong \mathtt C_2\times \mathtt C_2$.  The classifier is trained on a random sample of curves with conductor in the range specified by the first column, and verified on those which remain.
}}
\label{t:torsstructure}
\end{table}

Continuing with the torsion group, let us see how well the actual torsion group can be distinguished.
We established a balanced dataset of size $\sim 5 \times10^3~(\times2)$ for $\mathtt C_4$ and $\mathtt C_2 \times \mathtt C_2$ altogether.
Using a random forest classifier, we found that $E(\mathbb{Q})_{\mathrm{tors}}$ being $\mathtt C_4$ or $\mathtt C_2 \times \mathtt C_2$ can be separated using 500 Euler coefficients to fairly good accuracy.
The results are summarized in Table~\ref{t:torsstructure}. Note that the size of the dataset is relatively small compared to those of previous experiments. With a larger dataset, the precision might be improved. 

\subsection{Integral points}\label{sec:integral}

\begin{table}[h!!!]
\begin{center}
{\small \begin{tabular}{|c|c|c|c|c|}
\hline
$Q_E$ range& |Data| =  $\#\{E\}$ & $N$ &Precision&Confidence\\
\hline
$[1,5\times10^4]$&  $3.2\times10^4$ ($\times2$)& 500 &0.999&0.998\\
\hline
\end{tabular}
}
\end{center}
\caption{{\sf
The above table shows the precision and confidence of a naive Bayes classifier when asked to distinguish elliptic curves over $\mathbb{Q}$ with no integral points from those with a single integral point.  The classifier is trained on a random sample of curves with conductor in the range specified by the first column, and verified on those which remain.
}}
\label{t:integerpoints}
\end{table}

It is known that an elliptic curve has only finitely many integral points \cite[Chapter~VIII,~Chapter~IX,~Theorem~3.1]{silverman}. In contrast, it may have infinitely many rational points (this is the case when the rank $r_E>0$), as addressed above.
We set up a supervised ML to try to distinguish curves with no integral points from those with a single integral point, a total of around 60 thousand curves with conductor in the interval $[1,5 \times 10^4]$.
A balanced data-set of size $\sim 3.2 \times10^4~(\times2)$ for ``single integral point'' or ``no integral point'' was established and
a naive Bayes classifier produced the results summarized in Table~\ref{t:integerpoints}.
One can see that the results are extremely good. We will revisit this experiment in Section~\ref{sec:conjectures}.

\subsection{Tate--Shafarevich group}\label{sec:ECTS}

\begin{table}[h!!!]
\begin{center}
{\small \begin{tabular}{|c|c|c|c|c|}
\hline
$Q_E$ range&|Data| =  $\#\{E\}$ &$N$ & Precision\\
\hline
$[1,10^6]$ & $2.8\times10^4$ ($\times2$)&500
&<0.6\\
\hline
\end{tabular}
}
\end{center}
\caption{{\sf
The above table shows the precision of all classifiers when asked to distinguish elliptic curves over $\mathbb{Q}$ with Tate--Shafarevich order $4$ from those with order $9$.  The classifiers are trained on a random sample of curves with conductor in the range specified by the first column, and verified on those which remain. 
}}
\label{t:TSorder}
\end{table}

Finally, we come to the Tate--Shafarevich group, one of the most subtle parts of BSD. A definition of the Tate--Shafarevich group is given in \cite[Section~X.4]{silverman}. The Tate--Shafarevich group is not known to be finite, and there are no effective methods available for its computation. The LMFDB records the analytic order of the Tate--Shafarevich group, that is the real number implied by equation~\eqref{eq.BSDF}, which is equal to the order conditionally on the BSD conjecture. If the Tate--Shafarevich group is finite, then its order is a square integer.

We could try the following binary classification problem: take 500 Euler coefficients and see whether one could distinguish between a Tate--Shafarevich group of order 4 versus 9.
We tried a variety of methods, such as Bayesian or logistic classifiers, as well as some forward-feeding neural-networks,
but none performed especially well.
This is in accordance with the difficulty in computing this group.
The results are summarized in Table~\ref{t:TSorder}. 

For this problem alone, we implemented Weierstrass coefficient training (as was done in \cite{Alessandretti:2019jbs}). This experimental variant did not do well with any of the standard classifiers or regressors, again yielding no better than $<0.6$ precision. Nevertheless, we briefly review this approach for completeness.
Every elliptic curve over $\mathbb{Q}$ has a unique \textsl{reduced} minimal Weierstrass equation of the form:
\begin{equation}\label{eq.EllRedMinWeierstrass}
\begin{split}
y^2+e_1xy+e_2y=x^3+e_3x^2+e_4x+e_5,\\ 
e_1,e_3\in\{0,1\},\ \ e_2\in\{-1,0,1\}, \ \ e_4,e_5\in\mathbb{Z}.
\end{split}
\end{equation}
Using the coefficients in \eqref{eq.EllRedMinWeierstrass}, we define the vector:
\begin{equation}\label{eq.WeierstrassVector}
v_W(E)=(e_1,e_2,e_3,e_4,e_5)\in\mathbb{Z}^5.
\end{equation}
Let $\mathcal{F}$ denote a finite set of elliptic curves, and, for each $E\in\mathcal{F}$, let $I(E)$ be an invariant of interest. For example, $\mathcal{F}$ could be the set of all elliptic curves over $\mathbb{Q}$ with conductor less than one million and, for $E\in\mathcal{F}$, the invariant $I(F)$ could be the rank of $E$. We introduce the following labeled dataset: 
\begin{equation}\label{eq.WeierstrassData}
\mathcal{D}_W=\{v_W(E)\rightarrow I(E):E\in\mathcal{F}\}.
\end{equation}
Such a labeled dataset was used in \cite{Alessandretti:2019jbs}. 

\subsection{Interpretation of naive Bayesian models}\label{sec:conjectures}

Of the experimental results above, two instances with strikingly high accuracies are: the order of torsion subgroups in $E(\mathbb{Q})$ (Section~\ref{sec:torsorder}), and, the existence of integral points on $E$ (Section~\ref{sec:integral}).  The naive Bayes classifier was found to be optimal in both cases. Below we explore possible explanations.

We first observe that these classification problems are related to one another. Indeed, it can be shown that\footnote{We are grateful to \'Alvaro Lozano-Robledo and Chris Wuthrich, who informed the authors that one can prove these statements using the Nagell--Lutz theorem and other facts about elliptic curves.}:
\begin{enumerate}
\item If $\#E(\mathbb{Z})=1$, then $\#E(\mathbb{Q})_{\mathrm{tors}}=2$. Furthermore, the unique integral point is the torsion generator. 
\item If $\#E(\mathbb{Z})=0$, then $\#E(\mathbb{Q}_{\mathrm{tors}})\leq2$. Furthermore, in the LMFDB data, $99.99\%$ of optimal elliptic curves over $\mathbb{Q}$ with $\#E(\mathbb{Z})=0$ have torsion order 1.
\end{enumerate}
We might therefore expect that if a classifier can distinguish between $\#E(\mathbb{Q})_{\mathrm{tors}}\in\{1,2\}$ then it can distinguish between $\#E(\mathbb{Z})\in\{0,1\}$. 

On the other hand, we observe the following ``human'' procedure for distinguishing between torsion order $1$ and $2$ using the vectors $v_L(E)$ as in equation~\eqref{eq.Lvector}. Recall from equation~\eqref{eq.ap} that $a_p=p+1-\#E(\mathbb{F}_p)$. When $p$ is an odd prime, it follows that $a_p$ is even if and only if $\#E(\mathbb{F}_p)$ is even. If $p$ is moreover a prime of good reduction, then a point of order 2 in $E(\mathbb{Q})$ maps to a point of order 2 mod $p$ and so $\#E(\mathbb{F}_p)$ is even. We conclude that if $\#E(\mathbb{Q})_{\mathrm{tors}}=2$, then the vector $v_L(E)$ consists of even integers with a few possible exceptions coming from $p=2$ and bad primes (the exceptions are actually $\pm1$). In the case $\#E(\mathbb{Q})_{\mathrm{tors}}=1$ we observe that $a_p$'s are frequently odd as well as even. We are led to speculate that a naive Bayes classifier successfully distinguishes between vectors whose entries are all even from those whose entries are a mixture of even and odd numbers. 

To test this, we perform the following experiment. We generate one set of $100$-dimensional vectors with random integer coordinates in the range $[-10,10]$, and another set of $100$-dimensional vectors with coordinates equal to two times a random integer in the range $[-5,5]$. A naive Bayes classifier is able to distinguish these vectors to $99.8\%$ accuracy. By comparison, a random forest achieves around $74\%$. These accuracies are comparable to those observed in our experiments in Sections \ref{sec:torsorder} and \ref{sec:integral} and confirms the expectation that a Bayes classifier recognizes this difference.

\begin{figure}[h!!!]
(a)
$\begin{array}{c}\includegraphics[width=0.4\textwidth]{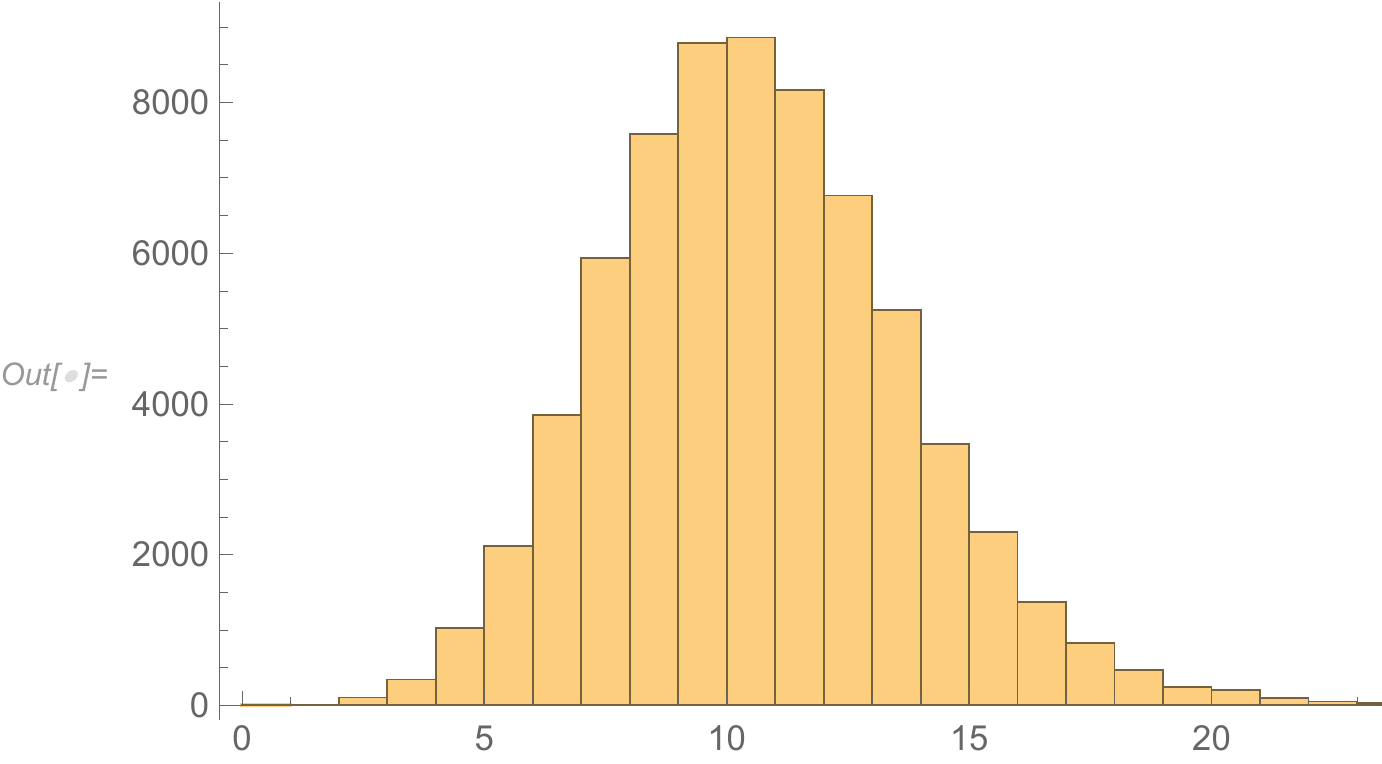}\end{array}$
(b)
$\begin{array}{c}\includegraphics[width=0.4\textwidth]{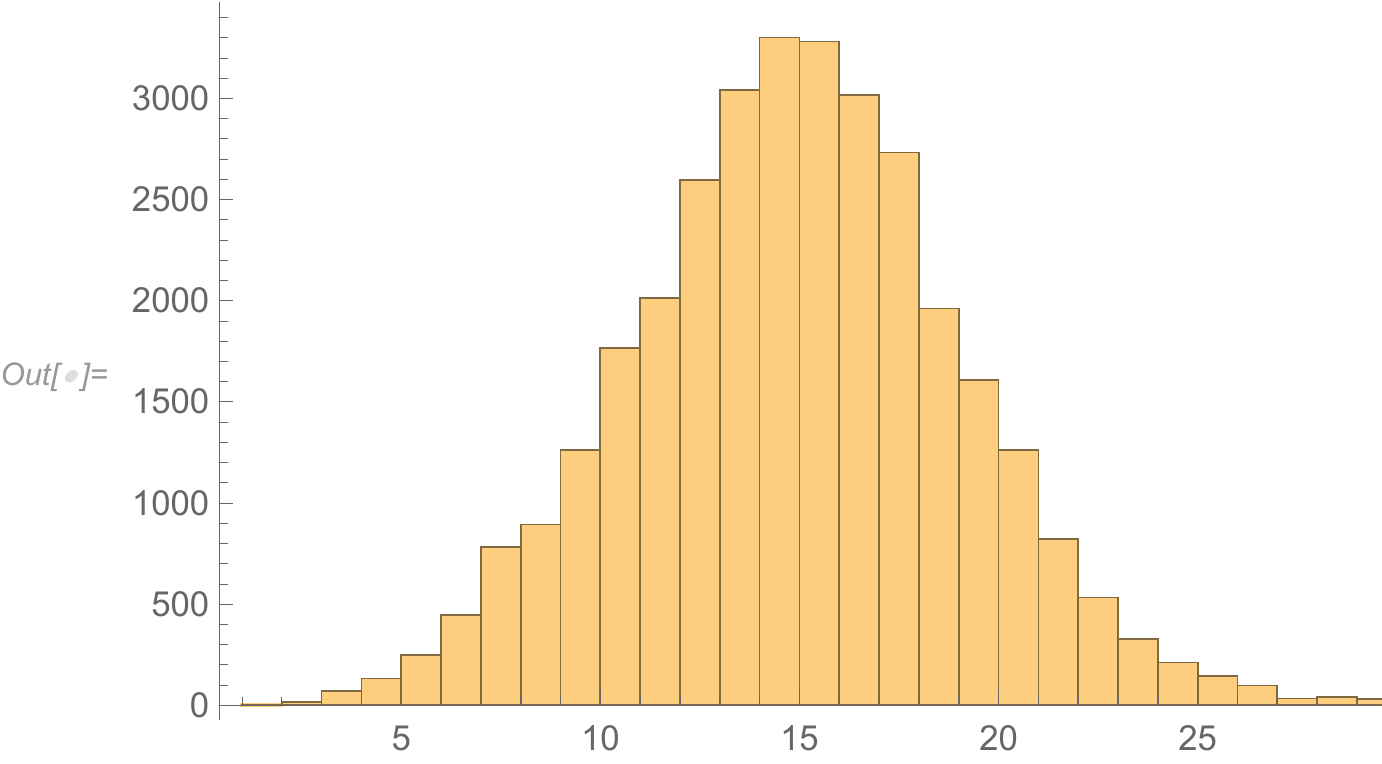}\end{array}$
\caption{
{\sf
The number of zeros in the list of 500 Euler coefficients for elliptic curves (a) without integral points, and (b) with a single integral point. 
}
}
\label{f:ap-0}
\end{figure}

As a separate attempt, inspired by the Lang--Trotter \cite{lt} conjecture, we study the distribution of the number of zeros in the vector $v_L$ of 500 $a_p$ values, for curves with a single integral point and without integral points (the case with the order of torsion group being 1 or 2 is similar to the ensuing discussions, mutatis mutandis).  For this purpose, we separately consider the set $\mathcal F_0$ of curves with no integral points and the set $\mathcal F_1$ of curves with a single integral point, where $|\mathcal F_0|=68,154$ and $|\mathcal F_1|=32,816$. We calculate the number of zeros in $v_L$ for each curve $E \in \mathcal F_i$, $i=0,1$, and draw the resulting histograms. This is shown in parts (a) and (b) respectively in Figure \ref{f:ap-0}. Clearly, the means of the two distributions are different. Precisely, $\mathcal F_0$ has mean $10.85$ with standard deviation $13.21$, while $\mathcal F_1$ has mean $15.26$ with standard deviation $14.70$.

To check whether a naive Bayes classifier detects this difference, we define the following binary vector for a positive integer $N$:
\begin{equation}\label{eq.binarytraining}
v_B(E)=(\delta_{p_1},\dots,\delta_{p_N})\in\{0,1\}^N,\ \ \delta_{p_i}=\begin{cases}0,&a_p=0,\\1,&a_p\neq0.\end{cases}
\end{equation}
The binary vectors in equation~\eqref{eq.binarytraining} are analogous to the binary vectors used in \cite{HLOb}. Replacing $v_L(E)$ with $v_B(E)$ in Section~\ref{sec:strategy} (Step 3) and performing the experiment in Section~\ref{sec:integral}, we observe that the naive Bayes is accurate to around $0.8$ precision. The result is similar if we instead use the following ternary vectors:
\begin{equation}\label{eq.ternarytraining}
v_T(E)=(\epsilon_{p_1},\dots,\epsilon_{p_N})\in\{-1,0,1\}^N,\ \ \epsilon_{p_i}=\begin{cases}-1,&a_p<0,\\0,&a_p=0,\\1,&a_p>0.\end{cases}
\end{equation}
Therefore, we see that what the Bayes classifier is picking up to reach the near 100\% predictions is based on {\it more than} merely the frequency of zeros/positives/negatives.
If we do include the actual values of the Euler coefficients, it takes at least around 7 coefficients to get to more than 0.9 accuracy.

On the other hand, we should point out that the number of zeros to the Euler coefficients is part of the Lang--Trotter \cite{lt} conjecture which is a refinement of the Sato--Tate \cite{Ta} conjecture. We are not aware of any claims in the literature that relate the distribution of zeros in the Euler coefficients to the number of integral points on or the torsion order of an elliptic curve. Our experimental results suggest that such relations may exist.
 
\section{Genus 2 curves}\label{s:genus2}

Having met with success for the genus 1 case,
in this section we describe our experimental results for genus $2$ curves defined over $\mathbb{Q}$. 
Throughout, we take the $N=200$ Euler coefficients and a conductor range from 1 to 1 million. This conductor range includes all the genus $2$ curves available in LMFDB. See \cite[Completeness~of~genus~2~curve~data~over~$\mathbb{Q}$]{lmfdb}. Again the $L$-function depends only on the isogeny class, but this time we cannot specify one curve per isogeny class in the LMFDB data.
Nevertheless, over 99\% of the genus 2 isogeny classes in the database (which has 65534 classes and 66158 curves) have a unique representative. Accepting this slight redundancy, we simply use all the genus 2 curves available in LMFDB.  

\subsection{Rank}
\begin{table}[h!!!]
\begin{center}
{\small \begin{tabular}{|c|c|c|}
\hline
\#\{$C$\}&Precision&Confidence\\
\hline
$1.21\times10^4$ ($\times 3$)&0.971&0.958\\
\hline
\end{tabular}
}
\end{center}
\caption{{\sf
The above table shows the precision and confidence of a logistic regression classifier when asked to distinguish between genus 2 curves over $\mathbb{Q}$ with rank in the set $\{0,1,2\}$.  
}}
\label{t:rank2}
\end{table}

We performed an experiment analogous to that in Section~\ref{sec:ECrank}. In the current context, a significant proportion of genus 2 curves have rank 2 and we consider the ternary classification problem of predicting whether the rank is 0, 1, or 2, from the Euler coefficients.
A balanced dataset of size $\sim 1 \times 10^4~(\times3)$ was thus established and a  logistic regression classifier was found to perform well, with accuracies $\sim 0.97$. We emphasize that this is a 3-way classification and to obtain this level of accuracies in impressive.
The results are summarized in Table~\ref{t:rank2}.

\subsection{Torsion order}

\begin{table}
\begin{center}
{\small \begin{tabular}{|c|c|c|}
\hline
$\#\{C\}$ &Precision&Confidence\\
\hline
                   $1.46\times10^4$ ($\times2$)&0.926&0.854\\
\hline
\end{tabular}
}
\end{center}
\caption{{\sf
The above table shows the precision and confidence of a naive Bayes classifier when asked to distinguish genus 2 curves over $\mathbb{Q}$ with torsion order $1$ from those with torsion order $2$.  
}}
\label{t:torsorder2}
\end{table}

As with the genus 1 case, we can try to distinguish the torsion group of order 1 versus 2 in a binary classification (cf. Section~\ref{sec:torsorder}).
A balanced data-set was established, with size $\sim 1.5 \times 10^4~(\times2)$ and a naive Bayes classifier was found to perform best, with results presented in Table~\ref{t:torsorder2}.

\subsection{Rational points}

\begin{table}
\begin{center}
{\small \begin{tabular}{|c|c|c|}
\hline
$\#$ of rational points &\#\{$C$\}&Precision\\
\hline
$\{0,1,2,3,4,5,6\}$& $5\times10^3~(\times 7) $&$<0.34$\\
\hline
$\{2,4\}$& $9.4\times10^3~(\times 2) $&$<0.73$\\
\hline
\end{tabular}
}
\end{center}
\caption{{\sf
The above table shows the precision of all classifiers when asked to distinguish between genus $2$ curves over $\mathbb{Q}$ with number of rational points as in the first column. 
}}
\label{t:rationalpoints}
\end{table}

As mentioned in the Introduction, curves of genus $>1$ have only a finite number of rational points. This allows for an experiment slightly different in nature to what was possible with elliptic curves.
Indeed, one could ask for a multi-category classification using the number of rational points, being predicted from the Euler coefficients. We tried various classes, after balancing the data but no classifier performed especially well. The results are summarized in Table~\ref{t:rationalpoints}, where a 7-way classification is shown in the first row, and a binary, in the second.
We suspect that training with a larger data set would result in a better performance. 

\subsection{Trivial Tate--Shafarevich group}

\begin{table}
\begin{center}
{\small \begin{tabular}{|c|c|c|}
\hline
\#\{$C$\}&Precision&Confidence\\
\hline
$4.2\times10^4$ ($\times 2$)&0.78&0.562\\
\hline
\end{tabular}
}
\end{center}
\caption{{\sf
The above table shows the precision and confidence of a logistic regression classifier when asked to distinguish genus 2 curves over $\mathbb{Q}$ with trivial Tate--Shafarevich group from those with non-trivial Tate--Shafarevich group.  
}}
\label{t:TS2}
\end{table}

Finally, we move to the  Tate--Shafarevich group.
Note that the order now needs not be a square for a genus 2 curve.
We performed a binary-classification (having established a balanced data set of size $\sim 4 \times10^4~(\times2)$) of whether Tate--Shafarevich group is trivial or not.
Again, no classifier was found to perform particularly well, though a logistic regression classifier performed best (see Table~\ref{t:TS2}), and the accuracies are comparable to those of the genus 1 case. Once again, the prediction is better than completely random. 

\section{Conclusions and Outlook}\label{s:outlook}

The experiments in this paper show that an ML classifier can be trained to predict the rank and the torsion order of an elliptic curve or a genus 2 curve with high precision when the curve is represented by a few hundred Euler coefficients. In particular, for elliptic curves, the torsion order and the number of integral points are determined almost perfectly by ML classifiers. Among the discrete invariants appearing in the BSD conjecture, only the order of the Tate--Shafarevich group seems to be out of reach with our approach of using a finite number of Euler coefficients.  

Along with our previous work \cite{HLOa, HLOb}, this paper confirms that ML classifiers perform surprisingly well with various invariants in number theory. High accuracies attained in our experiments reflect that data sets arising from mathematics are actually ``clean and clear'' without any noise. Prospectively, this opens up new opportunities of developing ML techniques for mathematics which exploit mathematical structures in data sets.

With all these experimental results and evidence at hand, a compelling call to action is to understand what ML classifiers actually recognize in the data sets. Though the algorithms of standard classifiers are well-known, it does not seem straightforward to precisely analyse what a classifier does with data sets. 

In another direction, we are reminded that the influential Langlands program anticipates correspondences between two kinds of data sets: {\em arithmetic} data and {\em automorphic} data. We have been experimenting with arithmetic data. In accordance with Langlands program, we expect that a machine would learn automorphic data with high precision and efficiency. It would be very interesting to investigate whether this expectation is valid. 


{\small 
Yang-Hui He {\sf hey@maths.ox.ac.uk} \\
Department of Mathematics, City, University of London, EC1V 0HB, UK;\\
Merton College, University of Oxford, OX14JD, UK;\\
School of Physics, NanKai University, Tianjin, 300071, P.R.~China

Kyu-Hwan Lee {\sf khlee@math.uconn.edu} \\
Department of Mathematics, University of Connecticut, Storrs, CT, 06269-1009, USA

Thomas Oliver {\sf Thomas.Oliver@nottingham.ac.uk} \\
School of Mathematical Sciences, University of Nottingham, University Park, \\
Nottingham, NG7 2QL, UK
}

\end{document}